\documentclass[a4paper,english,oneside]{article}
\usepackage[T1]{fontenc}
\usepackage[utf8]{inputenc} 

\usepackage{geometry}
\usepackage{amsmath}
\usepackage{amsfonts}
\usepackage{amssymb}
\usepackage{amsthm}
\usepackage[english]{babel}
\usepackage{color}
\usepackage{lmodern}
\usepackage{pstricks-add}
\usepackage{hyperref}
\hypersetup{hidelinks}

\usepackage{mathscinet}
\usepackage{enumitem}
\usepackage[cal=boondoxo,bb=ams]{mathalfa}
\usepackage{microtype}

\usepackage{mathtools}
\usepackage{esint}

\theoremstyle{plain}
\newtheorem{theorem}{Theorem}[section]
\newtheorem{lemma}[theorem]{Lemma}

\newtheorem{proposition}[theorem]{Proposition}

\theoremstyle{definition}

\theoremstyle{remark}
\newtheorem{remark}{Remark}

    \DeclareMathOperator\tr{Tr}
     
\begin{document}

\title{A global existence result for a semilinear wave equation with lower order terms on compact Lie groups}

\author{ Alessandro Palmieri} 

\date{}

\maketitle

\begin{abstract}

In this paper, we study the semilinear wave equation with lower order terms (damping and mass) and with power type nonlinearity $|u|^p$ on compact Lie groups. We will prove the  global in time existence of small data solutions in the evolution energy space without requiring any lower bounds for $p>1$. In our approach, we employ some results from Fourier analysis on compact Lie groups.
\end{abstract}

\begin{flushleft}
\end{flushleft}

\begin{flushleft}
\textbf{AMS Classification (2020)} Primary:  
35B45, 35L71; Secondary: 43A30, 43A77, 58J45
\end{flushleft}

\section{Introduction}

Let $\mathbb{G}$ be a compact, connected Lie group and let $\mathcal{L}$ be the Laplace -- Beltrami operator on $\mathbb{G}$ (which coincides with the Casimir element of the enveloping algebra). 
In the present work, we prove the global existence of small data solutions for the Cauchy problem for the semilinear wave equation with damping and mass and with power type nonlinearity, namely,
\begin{align}\label{semilinear wave CP}
\begin{cases} \partial_t^2 u-\mathcal{L} u+b\, \partial_t u+m^2u=| u|^p, &  x\in \mathbb{G}, \ t>0,\\
u(0,x)= u_0(x), & x\in \mathbb{G}, \\ \partial_t u(0,x)= u_1(x), & x\in \mathbb{G},
\end{cases}
\end{align} where $p>1$ and $b,m^2$ are positive constants.  

In \cite{Pal20WE,Pal20WEd} the semilinear wave equation in compact Lie groups ($b=m^2=0$) and the semilinear damped wave equation in compact Lie groups ($b>0,m^2=0$) have been considered, respectively. For both these cases, the nonexistence of global in time solutions can be proved for any $p>1$, under suitable sign assumptions on the Cauchy data. Purpose of this paper is to show that the combined presence of both a damping term $b\partial_t u$ and of a mass term $m^2 u$ modifies completely the situation. Indeed, we will show that, as in the Euclidean case (see \cite[Chapter 4]{PizMT} and \cite{ER18book}) or, more generally, as for a graded Lie group (cf. \cite{RT18}), the solution to the  linear homogeneous Cauchy problem
\begin{align}\label{linear wave CP}
\begin{cases} \partial_t^2 u-\mathcal{L} u+b\, \partial_t u+m^2u=0, &  x\in \mathbb{G}, \ t>0,\\
u(0,x)= u_0(x), & x\in \mathbb{G}, \\ \partial_t u(0,x)= u_1(x), & x\in \mathbb{G},
\end{cases}
\end{align} and its derivatives fulfill $L^2(\mathbb{G})$ -- $L^2(\mathbb{G})$ decay estimates with exponential decay rates with respect to the time variable. Using these decay estimates on $L^2(\mathbb{G})$ basis, we will be able to prove the existence of a uniquely determined globally in time defined mild solution to \eqref{semilinear wave CP} provided that $(u_0,u_1)$  has sufficiently small norm in the energy space via a standard contraction argument. Furthermore, a Gagliardo -- Nirenberg type inequality recently proved in \cite{RY19} will be employed to estimate the power nonlinearity in $L^2(\mathbb{G})$.

Let us stress that the decay estimates for the linear homogeneous Cauchy problem \eqref{linear wave CP} will be proved via Fourier Analysis on compact Lie groups. In particular, the group Fourier transform and Plancherel formula with respect to the space variables will play a crucial role.

Throughout the paper we denote by $L^q(\mathbb{G})$ the space of $q$ -- summable functions with respect to the normalized Haar measure on $\mathbb{G}$ for $1\leqslant q < \infty$ (respectively, essentially bounded for $q=\infty$) and for $s>0$ and $q\in(1,\infty)$ the Sobolev space $H^{s,q}_\mathcal{L}(\mathbb{G)}$ is defined as the space $$H^{s,q}_\mathcal{L}(\mathbb{G)}\doteq \left\{ f\in L^q(\mathbb{G}): (-\mathcal{L})^{s/2}f\in L^q(\mathbb{G})\right\}$$ equipped with the norm
$\|f\|_{H^{s,q}_\mathcal{L}(\mathbb{G)}}\doteq\| f\|_{L^q(\mathbb{G})}+\| (-\mathcal{L})^{s/2} f\|_{L^q(\mathbb{G})}.$ As customary, the Hilbert space $H^{s,2}_\mathcal{L}(\mathbb{G)}$ is simply denoted by $H^{s}_\mathcal{L}(\mathbb{G)}$.

Let us state the main result of this work, which is the global in time existence of small data solutions for  the semilinear Cauchy problem \eqref{semilinear wave CP}.

\begin{theorem}  \label{Thm glob esistence}
Let $\mathbb{G}$ be a compact, connected Lie group with  topological dimension $n$ such that $n\geqslant 3$. Let $(u_0,u_1)\in H_\mathcal{L}^1(\mathbb{G})\times L^2(\mathbb{G})$ and let $p>1$ such that $p\leqslant \frac{n}{n-2}$. Then, there exists $\varepsilon_0>0$ such that for any $\|(u_0,u_1)\|_{H_\mathcal{L}^1(\mathbb{G})\times L^2(\mathbb{G})}\leqslant \varepsilon_0$ the semilinear Cauchy problem \eqref{semilinear wave CP} admits a uniquely determined mild solution $$u\in \mathcal{C}\left([0,\infty),H^1_\mathcal{L}(\mathbb{G})\right)\cap \mathcal{C}^1\left([0,\infty),L^2(\mathbb{G})\right).$$ Furthermore, $u$ satisfies the decay estimates
\begin{equation} \label{decay estimate thm semilin}
\begin{split}
\|u(t,\cdot)\|_{L^2(\mathbb{G})} &\leqslant C d_{b,m^2}(t) \,  \|(u_0,u_1)\|_{H_\mathcal{L}^1(\mathbb{G})\times L^2(\mathbb{G})},  \\
\|(-\mathcal{L})^{1/2}u(t,\cdot)\|_{L^2(\mathbb{G})} &\leqslant C d_{b,m^2}(t) \, \|(u_0,u_1)\|_{H_\mathcal{L}^1(\mathbb{G})\times L^2(\mathbb{G})},  \\
\| \partial_t u(t,\cdot)\|_{L^2(\mathbb{G})} &\leqslant C d_{b,m^2}(t) \, \|(u_0,u_1)\|_{H_\mathcal{L}^1(\mathbb{G})\times L^2(\mathbb{G})}, 
\end{split}
\end{equation}
for any $t\geqslant 0$, where $C$ is a positive constant and the decay function $d_{b,m^2}(t)$ is given by 
\begin{align} \label{def decay functions}
d_{b,m^2}(t) \doteq \begin{cases} \mathrm{e}^{-\frac b2  t}  & \mbox{if} \ \  b^2< 4m^2, \\
 t\, \mathrm{e}^{-\frac b2  t}   & \mbox{if} \ \  b^2= 4m^2, \\
\mathrm{e}^{\big(-\frac b2 +\sqrt{\frac{b^2}{4}-m^2} \big) t}   & \mbox{if} \ \ b^2> 4m^2. \\ \end{cases}
\end{align}  

\end{theorem}

\begin{remark} \label{Remark ub p} The upper bound for the exponent $p$ in Theorem \ref{Thm glob esistence} is required in order to apply a Gagliardo -- Nirenberg type inequality  proved in \cite[Remark 1.7]{RY19}. The restriction $n\geqslant 3$ on the topological dimension is made to fulfill the conditions for the application of this inequality too and it can be removed by looking for solutions in spaces with less regularity, namely, in $\mathcal{C}\left([0,\infty),H^s_\mathcal{L}(\mathbb{G})\right)\cap \mathcal{C}^1\left([0,\infty),L^2(\mathbb{G})\right)$ for some $s\in (0,1)$.
\end{remark}

\begin{remark} Let us stress that the decay estimates for the solution of the semilinear Cauchy problem are exactly the same ones as for the solution of the corresponding linear homogeneous problem (cf. Proposition \ref{Prop L^2-L^2 estimates}). This is a consequence of the construction of the family of evolution spaces $\{X(T)\}_{T>0}$ equipped with a suitable weighted norm in the proof of Theorem \ref{Thm glob esistence} (cf. \eqref{def norm X(T)} in Section \ref{Section approach}).
\end{remark}

\begin{remark} 
Note that in the statement of Theorem \ref{Thm glob esistence} no further lower bound for $p>1$ is required. This remarkable property is due to the exponential decay rates in the $L^2(\mathbb{G})$ -- $L^2(\mathbb{G})$ estimates for the solution to the corresponding linear homogeneous Cauchy problem from Proposition \ref{Prop L^2-L^2 estimates} below.
\end{remark}

\subsection*{Notations} Throughout this paper the following notations are used: as in the introduction, $\mathcal{L}$ denotes the Laplace -- Beltrami operator on $\mathbb{G}$; $\tr(A)= \sum_{j=1}^d a_{jj}$ denotes the trace of the matrix $A=(a_{ij})_{1\leqslant i,j\leqslant d} \in \mathbb{C}^{d\times d}$ while $A^*=(\overline{a_{ji}})_{1\leqslant i,j\leqslant d}$ denotes its adjoint matrix; $I_d\in\mathbb{C}^{d\times d}$ is the identity matrix; $\mathrm{d}x$ stands for the normalized Haar measure on the compact group $\mathbb{G}$; finally, we write $f\lesssim g$ when there exists a positive constant $C$ such that $f\leqslant Cg$. 

\section{Philosophy of our approach} \label{Section approach}

In this section we will clarify our strategy to prove Theorem \ref{Thm glob esistence}. First, we provide the notion of mild solutions to \eqref{semilinear wave CP} in our framework. By applying Duhamel's principle, the solution to the linear problem
\begin{align}\label{linear CP wave inhom}
\begin{cases} \partial_t^2 u-\mathcal{L} u+b\, \partial_t u+m^2u =F(t,x), &  x\in \mathbb{G}, \ t>0,\\
u(0,x)= u_0(x), & x\in \mathbb{G}, \\ \partial_t u(0,x)= u_1(x), & x\in \mathbb{G}
\end{cases}
\end{align}  can be represented as
\begin{align*}
u(t,x) = u_0(x)\ast_{(x)} E_0(t,x;b,m^2) + u_1(x)\ast_{(x)} E_1(t,x;b,m^2)+ \int_0^t F(s,x)\ast_{(x)} E_1(t-s,x;b,m^2) \, \mathrm{d}s.
\end{align*}
 where $E_0(t,x;b,m^2)$ and $E_1(t,x;b,m^2)$ denote, respectively, the fundamental solutions to \eqref{linear CP wave inhom} in the homogeneous case $F= 0$ with initial data $(u_0,u_1)=(\delta_0,0)$ and $(u_0,u_1)=(0,\delta_0)$.
Let us point out that in order to get the previous representation formula we applied the invariance by time translations for the differential  operator $\partial_t^2-\mathcal{L}+b\, \partial_t+m^2$ and the property $L\big(v\ast_{(x)}E_1(t,\cdot\, ;b,m^2)\big)=v\ast_{(x)}L(E_1(t,\cdot\, ;b,m^2))$ for any left -- invariant differential operator $L$ on $\mathbb{G}$.

The function $u$ is said a \emph{mild solution} to \eqref{semilinear wave CP} on $[0,T]$ if $u$ is a fixed point for the nonlinear integral operator $N$ defined by
\begin{align*}
N: u\in X(T) \to N u(t,x) & \doteq  u_0(x)\ast_{(x)} E_0(t,x;b,m^2)  +
 u_1(x)\ast_{(x)} E_1(t,x;b,m^2) \\ & \qquad+ \int_0^t |u(s,x)|^p\ast_{(x)} E_1(t-s,x;b,m^2) \, \mathrm{d}s
\end{align*} on the evolution space $X(T)\doteq\mathcal{C}\left([0,T],H^1_\mathcal{L}(\mathbb{G})\right)\cap \mathcal{C}^1\left([0,T],L^2(\mathbb{G})\right)$, endowed with the norm
\begin{align} \label{def norm X(T)}
\|u\|_{X(T)}\doteq \sup_{t\in[0,T]} \left(d_{b,m^2}(t)\right)^{-1}\left(\|u(t,\cdot)\|_{L^2(\mathbb{G})}+\|(-\mathcal{L})^{1/2}u(t,\cdot)\|_{L^2(\mathbb{G})}+\|\partial_t u(t,\cdot)\|_{L^2(\mathbb{G})}\right).
\end{align}

If we show the validity of the inequalities
\begin{align}
\| Nu\|_{X(T)} & \leqslant C \| (u_0,u_1)\|_{H^1_\mathcal{L}(\mathbb{G})\times L^2(\mathbb{G})} + C \| u\|_{X(T)}^p, \label{1st ineq N op} \\
\| Nu -Nv\|_{X(T)} & \leqslant   C\| u-v\|_{X(T)}\left( \| u\|_{X(T)}^{p-1}+\| v\|_{X(T)}^{p-1}\right), \label{2nd ineq N op}
\end{align} for any $u,v\in X(T)$ and for a suitable constant $C>0$ independent of $T$, then, by Banach's fixed point theorem it follows that $N$ admits a uniquely determined fixed point $u$ provided that $\| (u_0,u_1)\|_{H^1_\mathcal{L}(\mathbb{G})\times L^2(\mathbb{G})}$ is small enough. This function $u$ is our mild solution to \eqref{semilinear wave CP} on $[0,T]$. Furthermore, thanks to the fact that \eqref{1st ineq N op} and \eqref{2nd ineq N op} hold uniformly  with respect to $T$, this solution can be prolonged and defined for any $t\in(0,\infty)$.

 Before considering the semilinear Cauchy problem, we will determine $L^2(\mathbb{G})$ -- $L^2(\mathbb{G})$ estimates for the solution to  \eqref{linear wave CP} via the group Fourier transform with respect to the spatial variable. 
We point out explicitly that the definition of the norm on the space $X(T)$ in \eqref{def norm X(T)} is compatible with the decay estimates for the solution to the linear Cauchy problem \eqref{linear wave CP} that we will determine afterwards (see also \cite[Section 18.1.1]{ER18book} for a detailed overview on this topic). 
 After that these estimates will have been established, we could show the global existence of small data solutions for \eqref{semilinear wave CP}  by applying a Gagliardo -- Nirenberg type inequality derived recently in \cite{RY19} (cf. Lemma \ref{Lemma GN ineq L2} in Section \ref{Subsection Fixed Point}). 
 
 The remaining part of the paper is organized as follows: in Section \ref{Subsection GFT} we recall the main tools from Fourier Analysis and representation theory on compact Lie groups which are necessary for our approach; in Section \ref{Subsection L^2-L^2 est}, we establish $L^2(\mathbb{G})$ -- $L^2(\mathbb{G})$ estimates for the solution of \eqref{linear wave CP} and its first order derivatives; in Section \ref{Subsection Fixed Point} it will be shown that the operator $N$ is a contraction on $X(T)$ provided that the data are sufficiently small in the energy space; finally, we provide some concluding remarks in Section \ref{Section final rem}.

\section{Group Fourier transform} \label{Subsection GFT}

Let us recall some results on Fourier Analysis on compact Lie groups. For further details on this topic we refer to the monographs \cite{RT10,FR16}. 

A \emph{continuous unitary representation} $\xi:\mathbb{G}\to\mathbb{C}^{d_\xi\times d_\xi}$ of dimension $d_\xi$ is a continuous group homomorphism from $\mathbb{G}$ to the group of unitary matrix $\mathrm{U}(d_\xi,\mathbb{C})$, that is, $\xi(xy)=\xi(x)\xi(y)$ and $\xi(x)^*=\xi(x)^{-1}$  for all $x,y\in\mathbb{G}$ and the elements $\xi_{ij}:\mathbb{G}\to \mathbb{C}$ of the matrix representation $\xi$ are continuous functions for all $i,j\in \{1,\ldots,d_\xi\}$. Two representations $\xi,\eta$ of $\mathbb{G}$ are said \emph{equivalent} if there exists an invertible intertwining operator $A$ such that $A \xi(x) = \eta(x) A$ for any $x\in \mathbb{G}$. A subspace $W\subset \mathbb{C}^{d_\xi}$ is $\xi$ -- invariant if $\xi(x) \cdot W\subset W$ for any $x\in \mathbb{G}$. A representation $\xi$ is \emph{irreducible} if the only $\xi$ -- invariant subspaces are the trivial ones $\{0\},\mathbb{C}^{d_\xi}$.

The unitary dual $\widehat{\mathbb{G}}$ of the compact Lie group $\mathbb{G}$ consists of the equivalence class $[\xi]$ of continuous irreducible unitary representation $\xi:\mathbb{G}\to\mathbb{C}^{d_\xi\times d_\xi}$. 

Given $f\in L^1(\mathbb{G})$, its Fourier coefficients at $[\xi]\in\widehat{\mathbb{G}}$ is defined by
\begin{align*}
\widehat{f}(\xi)\doteq \int_{\mathbb{G}} f(x) \xi(x)^*  \mathrm{d}x \in \mathbb{C}^{d_\xi\times d_\xi},
\end{align*} where the integral is taken with respect to the Haar measure on $\mathbb{G}$.

If $f\in L^2(\mathbb{G})$, then, the Fourier series representation for $f$ is given by
$$f(x) = \sum_{[\xi]\in\widehat{\mathbb{G}}} d_\xi \tr\big(\xi(x)\widehat{f}\ \big),$$
 where hereafter just one irreducible unitary matrix representation is picked in the sum for each equivalence class $[\xi]$ in $\widehat{\mathbb{G}}$.
Furthermore, for $f\in L^2(\mathbb{G})$ Plancherel formula takes the following form
\begin{align} \label{Plancherel formula G}
\| f\|^2_{L^2(\mathbb{G})} =  \sum_{[\xi]\in\widehat{\mathbb{G}}} d_\xi \big\| \widehat{f}(\xi)\big\|^2_{\mathrm{HS}}, 
\end{align}  where the Hilbert -- Schmidt norm of the matrix $\widehat{f}(\xi)$ is defined as follows:
\begin{align*}
 \big\| \widehat{f}(\xi)\big\|^2_{\mathrm{HS}}  \doteq  \tr\big( \widehat{f}(\xi) \widehat{f}(\xi)^*\big) = \sum_{i,j=1}^{d_\xi} \big| \widehat{f}(\xi)_{ij}\big|^2.
\end{align*}

For our analysis it is important to understand the behavior of the group Fourier transform  with respect to the Laplace -- Beltrami operator $\mathcal{L}$. Given $[\xi]\in\widehat{\mathbb{G}}$, then, all $\xi_{ij}$ are eigenfunctions for $\mathcal{L}$ with the same not positive eigenvalue $-\lambda_\xi^2$, namely,
\begin{align*}
-\mathcal{L} \xi_{ij}(x) = \lambda^2_\xi \, \xi_{ij}(x) \quad\mbox{for any} \ x\in\mathbb{G} \ \mbox{and for all} \ i,j\in \{1,\ldots,d_\xi\}.
\end{align*} This means that the symbol of $\mathcal{L}$ is 
\begin{align} \label{symbol Laplace-Beltrami}
\sigma_\mathcal{L} (\xi)= - \lambda^2_\xi I_{d_\xi},
\end{align} that is, $\widehat{\mathcal{L} f}(\xi)= \sigma_\mathcal{L} (\xi) \widehat{f}(\xi) =- \lambda^2_\xi \widehat{f}(\xi)$ for any $[\xi]\in\widehat{\mathbb{G}}$.

Finally, through Plancherel formula for $s>0$ we have
\begin{align*}
\|f\|^2_{\dot{H}^s_\mathcal{L}(\mathbb{G})} =\| (-\mathcal{L})^{s/2}f\|^2_{L^2(\mathbb{G})} =  \sum_{[\xi]\in\widehat{\mathbb{G}}} d_\xi \lambda_\xi^{2s}\big\| \widehat{f}(\xi)\big\|^2_{\mathrm{HS}}.
\end{align*}

\section{$L^2(\mathbb{G})$ -- $L^2(\mathbb{G})$ estimates for the solution to the linear \\ homogeneous problem} \label{Subsection L^2-L^2 est}

In this section, we derive $L^2(\mathbb{G})$ -- $L^2(\mathbb{G})$ estimates for the solution to the linear Cauchy problem \eqref{linear wave CP}.
 We follow the approach from \cite{GR15}, which have been  recently applied to study other semilinear hyperbolic models in compact Lie groups in \cite{Pal20WEd,Pal20WE}. The main idea of this approach is to employ the group Fourier transform with respect to the spatial variable $x$ to get an explicit expression for the $L^2(\mathbb{G})$ norms of $u(t,\cdot)$, $(-\mathcal{L})^{1/2}u(t,\cdot)$ and $\partial_t u(t,\cdot)$, respectively. Plancherel formula in the framework of compact Lie groups plays a fundamental role in this step.

\begin{proposition} \label{Prop L^2-L^2 estimates}
Let $\mathbb{G}$ be a compact Lie group. Let us assume $(u_0,u_1)\in H^1_\mathcal{L}(\mathbb{G})\times L^2(\mathbb{G})$ and let $u\in \mathcal{C}\big([0,\infty),H^1_\mathcal{L}(\mathbb{G})\big)\cap  \mathcal{C}^1\big([0,\infty),L^2(\mathbb{G})\big)$ be the solution to the homogeneous Cauchy problem \eqref{linear wave CP}.
 Then, the following $L^2(\mathbb{G})$ -- $L^2(\mathbb{G})$  decay estimates are satisfied
\begin{align}
\|u(t,\cdot)\|_{L^2(\mathbb{G})} &\leqslant C d_{b,m^2}(t)  \Big( \|u_0\|_{L^2(\mathbb{G})}+   \|u_1\|_{L^2(\mathbb{G})}\Big), \label{L^2 norm u(t)} \\
\|(-\mathcal{L})^{1/2}u(t,\cdot)\|_{L^2(\mathbb{G})} &\leqslant C d_{b,m^2}(t) \left( \|u_0\|_{H^1_\mathcal{L}(\mathbb{G})}+\|u_1\|_{L^2(\mathbb{G})}\right), \label{L^2 norm (-L)^1/2 u(t)} \\
\| \partial_t u(t,\cdot)\|_{L^2(\mathbb{G})} &\leqslant C d_{b,m^2}(t) \left( \|u_0\|_{H^1_\mathcal{L}(\mathbb{G})}+\|u_1\|_{L^2(\mathbb{G})}\right), \label{L^2 norm ut(t)}
\end{align} for any $t\geqslant 0$, where $C$ is a positive multiplicative constant and the decay function $d_{b,m^2}$ is defined in \eqref{def decay functions}.
\end{proposition}

\begin{remark} In the threshold case $b^2=4m^2$ we me refine \eqref{L^2 norm (-L)^1/2 u(t)}, considering  $\mathrm{e}^{-\frac b2  t}$ as decay rate. However, for the sake of brevity we skipped this detail, since in the proof of the global existence result the decay rate in \eqref{L^2 norm (-L)^1/2 u(t)} suffices to apply Banach's fixed point theorem for any $p>1$ without requiring a further lower bound for $p$.
\end{remark}

\begin{remark} We point out that no connectedness is required for $\mathbb{G}$ in the study of the linear Cauchy problem \eqref{linear wave CP} in the statement of Proposition \ref{Prop L^2-L^2 estimates}.
\end{remark}

\subsubsection*{Representation formula for the Fourier coefficient $\widehat{u}(t,\xi)_{k\ell}$}
Let $u$ solve \eqref{linear wave CP}. By $\widehat{u}(t,\xi) = (\widehat{u}(t,\xi)_{k\ell})_{1\leqslant k,\ell\leqslant d_\xi}\in \mathbb{C}^{d_\xi\times d_\xi}$, $[\xi]\in\widehat{\mathbb{G}}$ we denote the group Fourier transform of $u$ with respect to the $x$ -- variable. Therefore, $\widehat{u}(t,\xi)$ is a solution of the Cauchy problem for the system of ODEs (with the size of the system that depends on the representation $\xi$)
\begin{align*}
\begin{cases}
\partial_t^2 \widehat{u}(t,\xi) -\sigma_\mathcal{L}(\xi) \widehat{u}(t,\xi)+b\, \partial_t \widehat{u}(t,\xi) +m^2\, \widehat{u}(t,\xi) =0, & t>0, \\
 \widehat{u}(0,\xi) = \widehat{u}_0(\xi), \\
 \partial_t \widehat{u}(0,\xi) =  \widehat{u}_1(\xi).
\end{cases}
\end{align*} By using the symbol of the Laplace --  Beltrami operator in \eqref{symbol Laplace-Beltrami}, we obtain that the previous system is actually decoupled in $d_\xi^2$ independent scalar ODEs, namely,
\begin{align}\label{scalar dec ODE}
\begin{cases}
\partial_t^2 \widehat{u}(t,\xi)_{k\ell} + \lambda_\xi^2 \widehat{u}(t,\xi)_{k\ell} +b\, \partial_t \widehat{u}(t,\xi)_{k\ell}+m^2\, \widehat{u}(t,\xi)_{k\ell}=0, & t>0, \\
 \widehat{u}(0,\xi)_{k\ell} = \widehat{u}_0(\xi)_{k\ell}, \\
 \partial_t \widehat{u}(0,\xi)_{k\ell} =  \widehat{u}_1(\xi)_{k\ell},
\end{cases}
\end{align} for any $k,\ell\in \{1,\ldots,d_\xi\}$.
Straightforward computations lead to the representation formula 
\begin{align}\label{representation u hat kl}
  \widehat{u}(t,\xi)_{k\ell} =  \mathrm{e}^{-\frac{b}{2}t} G_0(t;b,m^2;\xi) \, \widehat{u}_0(\xi)_{k\ell} + \mathrm{e}^{-\frac{b}{2}t}G_1(t;b,m^2;\xi)\Big( \widehat{u}_1(\xi)_{k\ell}+\tfrac b2  \widehat{u}_0(\xi)_{k\ell}\Big)
\end{align} for the solution to the linear homogeneous Cauchy problem \eqref{scalar dec ODE},  where
\begin{equation} \label{def G0 G1}
\begin{split}
 G_0(t;b,m^2;\xi)  & \doteq \begin{cases} \cosh\left(  \sqrt{\tfrac{b^2}{4}-m^2- \lambda_\xi^2} \, t\right)  & \mbox{if} \ \lambda_\xi^2<\tfrac{b^2}{4}-m^2 , \\
 1  & \mbox{if} \ \lambda_\xi^2=\tfrac{b^2}{4}-m^2, \\
 \cos\left( \sqrt{\lambda_\xi^2+m^2-\tfrac{b^2}{4}} \, t\right)  & \mbox{if} \ \lambda_\xi^2>\tfrac{b^2}{4}-m^2,   \end{cases} \\
 G_1(t;b,m^2;\xi) & \doteq \begin{cases} \dfrac{\sinh\left( \sqrt{\tfrac{b^2}{4}-m^2- \lambda_\xi^2} \, t\right)}{\sqrt{\tfrac{b^2}{4}-m^2- \lambda_\xi^2}}  & \mbox{if} \ \lambda_\xi^2<\tfrac{b^2}{4}-m^2 , \\
 t  & \mbox{if} \ \lambda_\xi^2=\tfrac{b^2}{4}-m^2, \\
 \dfrac{\sin\left( \sqrt{\lambda_\xi^2+m^2-\tfrac{b^2}{4}} \, t\right)}{\sqrt{\lambda_\xi^2+m^2-\tfrac{b^2}{4}}}   & \mbox{if} \ \lambda_\xi^2>\tfrac{b^2}{4}-m^2.   \end{cases}
 \end{split}
\end{equation} Notice that $G_0(t;b,m^2;\xi)=\partial_t G_1(t;b,m^2;\xi)$ for any $[\xi]\in\widehat{\mathbb{G}}$.

Finally, by \eqref{representation u hat kl} for any $[\xi]\in\widehat{\mathbb{G}}$ and any $k,\ell\in\{1,\ldots,d_\xi\}$ we derive the following representation for the time derivative of $\widehat{u}(t,\xi)_{k\ell}$
\begin{align} \label{representation u_t hat kl}
\partial_t \widehat{u}(t,\xi)_{k\ell} =  \mathrm{e}^{-\frac{b}{2}t} G_0(t;b,m^2;\xi) \, \widehat{u}_1(\xi)_{k\ell} - \mathrm{e}^{-\frac{b}{2}t}G_1(t;b,m^2;\xi)\Big( \tfrac b2 \, \widehat{u}_1(\xi)_{k\ell}+ (\lambda_\xi^2+m^2)  \, \widehat{u}_0(\xi)_{k\ell}\Big).
\end{align} 

Let us prove now Proposition \ref{Prop L^2-L^2 estimates}. We will consider separately the three subcases $b^2 \lesseqqgtr 4m^2$.

\subsubsection*{Case $b^2<4m^2$}

 In this case the characteristic roots are always complex conjugate and in the representation formula \eqref{representation u hat kl} it has sense to consider only the case $\lambda_\xi^2 > \frac{b^2}{4}-m^2$, since all eigenvalues $\{\lambda_\xi^2\}_{[\xi]\in\widehat{\mathbb{G}}}$ of $-\mathcal{L}$ are nonnegative. So, we can estimate
\begin{equation} \label{estimate u hat kl case b^2<4m^2}
\begin{split}
|\widehat{u}(t,\xi)_{k\ell}| & \lesssim \mathrm{e}^{-\frac b2 t} \big(|\widehat{u}_0(\xi)_{k\ell}| +  |\widehat{u}_1(\xi)_{k\ell}|\big), \\
\lambda_\xi|\widehat{u}(t,\xi)_{k\ell}| & \lesssim \mathrm{e}^{-\frac b2 t} \big((1+	\lambda_\xi)|\widehat{u}_0(\xi)_{k\ell}| +  |\widehat{u}_1(\xi)_{k\ell}|\big), 
\end{split}
\end{equation} for any $t\geqslant 0$. Similarly, by \eqref{representation u_t hat kl} we get 
\begin{align}  \label{estimate u_t hat kl case b^2<4m^2}
|\partial_t \widehat{u}(t,\xi)_{k\ell}| & \lesssim \mathrm{e}^{-\frac b2 t} \big((1+	\lambda_\xi)|\widehat{u}_0(\xi)_{k\ell}| +  |\widehat{u}_1(\xi)_{k\ell}|\big)
\end{align}  for any $t\geqslant 0$.
Combining \eqref{estimate u hat kl case b^2<4m^2} and \eqref{estimate u_t hat kl case b^2<4m^2} and using Plancherel formula, we have
\begin{align}
\| (-\mathcal{L})^{i/2}\partial_t^j u(t,\cdot)\|^2_{L^2(\mathbb{G)}} & = \sum_{[\xi]\in\widehat{\mathbb{G}}} d_\xi \sum_{k,\ell=1}^{d_\xi} \lambda_\xi^{2i}|\partial_t^j\widehat{u}(t,\xi)_{k\ell}|^2 \notag \\ & \lesssim  \mathrm{e}^{- b t}  \sum_{[\xi]\in\widehat{\mathbb{G}}} d_\xi \sum_{k,\ell=1}^{d_\xi} \Big( \big(1+\lambda_\xi^{2}\big)^{(i+j)}|\widehat{u}_0(\xi)_{k\ell}|^2 +|\widehat{u}_1(\xi)_{k\ell}|^2 \Big) \notag \\ &  =  \mathrm{e}^{- b t} \Big( \| u_0\|^2_{H^{i+j}_\mathcal{L}(\mathbb{G)}} +\| u_1\|^2_{L^2(\mathbb{G)}}\Big) \label{proof L^2 est u(t) case b^2<4m^2}
\end{align} for any $i,j\in\{0,1\}$ such that $0\leqslant i+j\leqslant 1$ and any $t\geqslant 0$. Clearly \eqref{proof L^2 est u(t) case b^2<4m^2} implies \eqref{L^2 norm u(t)}, \eqref{L^2 norm (-L)^1/2 u(t)} and \eqref{L^2 norm ut(t)} in the case $b^2<4m^2$ (here we used $H^{0}_\mathcal{L}(\mathbb{G)}=L^2(\mathbb{G})$).

\subsubsection*{Case $b^2=4m^2$}

In this case the representation formula \eqref{representation u hat kl} is just
\begin{align*}
 \widehat{u}(t,\xi)_{k\ell} & =  \mathrm{e}^{-\frac{b}{2}t} \cos(\lambda_\xi t) \, \widehat{u}_0(\xi)_{k\ell} + \mathrm{e}^{-\frac{b}{2}t}\frac{\sin(\lambda_\xi t)}{\lambda_\xi}\Big( \widehat{u}_1(\xi)_{k\ell}+\tfrac b2  \widehat{u}_0(\xi)_{k\ell}\Big) & \mbox{for} \ \ [\xi]\in\widehat{\mathbb{G}} \ \ \mbox{such that} \ \ \lambda_\xi^2>0, \\
 \widehat{u}(t,\xi)_{k\ell} & =  \mathrm{e}^{-\frac{b}{2}t} \, \widehat{u}_0(\xi)_{k\ell} + t\, \mathrm{e}^{-\frac{b}{2}t}\Big( \widehat{u}_1(\xi)_{k\ell}+\tfrac b2  \widehat{u}_0(\xi)_{k\ell}\Big) & \mbox{for} \ \ [\xi]\in\widehat{\mathbb{G}} \ \ \mbox{such that} \ \ \lambda_\xi^2=0.
\end{align*} Note that we have to consider necessarily the second case, since for the trivial 1 -- dimensional representation $1:x\in\mathbb{G}\to 1\in\mathbb{C}$ we have $-\mathcal{L} (1)=0$. Therefore,  
\begin{equation} \label{estimate u hat kl case b^2=4m^2}
\begin{split}
|\widehat{u}(t,\xi)_{k\ell}| & \lesssim t\, \mathrm{e}^{-\frac b2 t} \big(|\widehat{u}_0(\xi)_{k\ell}| +  |\widehat{u}_1(\xi)_{k\ell}|\big), \\
\lambda_\xi|\widehat{u}(t,\xi)_{k\ell}| & \lesssim \mathrm{e}^{-\frac b2 t} \big((1+	\lambda_\xi)|\widehat{u}_0(\xi)_{k\ell}| +  |\widehat{u}_1(\xi)_{k\ell}|\big), 
\end{split}
\end{equation} for any $t\geqslant 0$. In the second inequality, we can drop the factor $t$ as the spectrum of $-\mathcal{L}$ is discrete and has no finite cluster points.

Similarly, from \eqref{representation u_t hat kl} it follows
\begin{align}  \label{estimate u_t hat kl case b^2=4m^2}
|\partial_t \widehat{u}(t,\xi)_{k\ell}| & \lesssim t \, \mathrm{e}^{-\frac b2 t} \big((1+	\lambda_\xi)|\widehat{u}_0(\xi)_{k\ell}| +  |\widehat{u}_1(\xi)_{k\ell}|\big)
\end{align}  for any $t\geqslant 0$.  Applying Plancherel formula twice as in \eqref{proof L^2 est u(t) case b^2<4m^2} and using \eqref{estimate u hat kl case b^2=4m^2} and \eqref{estimate u_t hat kl case b^2=4m^2} to control the Fourier coefficients, we obtain \eqref{L^2 norm u(t)}, \eqref{L^2 norm (-L)^1/2 u(t)} and \eqref{L^2 norm ut(t)} in the case $b^2=4m^2$.

\subsubsection*{Case $b^2>4m^2$}

In this last case, the characteristic roots may be either complex conjugate or coincident or real distinct, depending on the range for $\lambda_\xi^2$. Comparing all possible cases in \eqref{def G0 G1}, we get 
\begin{equation} \label{estimate u hat kl case b^2>4m^2}
\begin{split}
|\widehat{u}(t,\xi)_{k\ell}| & \lesssim  \mathrm{e}^{\big(-\frac b2+\sqrt{\frac{b^2}{4}-m^2}\big) t} \big(|\widehat{u}_0(\xi)_{k\ell}| +  |\widehat{u}_1(\xi)_{k\ell}|\big), \\
\lambda_\xi|\widehat{u}(t,\xi)_{k\ell}| & \lesssim \mathrm{e}^{\big(-\frac b2+\sqrt{\frac{b^2}{4}-m^2}\big) t} \big((1+	\lambda_\xi)|\widehat{u}_0(\xi)_{k\ell}| +  |\widehat{u}_1(\xi)_{k\ell}|\big), 
\end{split}
\end{equation} for any $t\geqslant 0$. We point out that the regularity is provided from the case with complex conjugate characteristic roots whereas the decay rates is given by the continuous irreducible unitary representations with $\lambda_\xi^2=0$.

Analogously, from \eqref{representation u_t hat kl} we get 
\begin{align}  \label{estimate u_t hat kl case b^2>4m^2}
|\partial_t \widehat{u}(t,\xi)_{k\ell}| & \lesssim \mathrm{e}^{\big(-\frac b2+\sqrt{\frac{b^2}{4}-m^2}\big) t} \big((1+	\lambda_\xi)|\widehat{u}_0(\xi)_{k\ell}| +  |\widehat{u}_1(\xi)_{k\ell}|\big)
\end{align}  for any $t\geqslant 0$.

Consequently, by Plancherel formula combined with \eqref{estimate u hat kl case b^2>4m^2} and \eqref{estimate u_t hat kl case b^2>4m^2}, we find
\begin{align}
\| (-\mathcal{L})^{i/2}\partial_t^j u(t,\cdot)\|^2_{L^2(\mathbb{G)}} & = \sum_{[\xi]\in\widehat{\mathbb{G}}} d_\xi \sum_{k,\ell=1}^{d_\xi} \lambda_\xi^{2i}|\partial_t^j\widehat{u}(t,\xi)_{k\ell}|^2 \notag \\ & \lesssim  \mathrm{e}^{\left(- b+\sqrt{b^2-4m^2}\right) t}  \sum_{[\xi]\in\widehat{\mathbb{G}}} d_\xi \sum_{k,\ell=1}^{d_\xi} \Big( \big(1+\lambda_\xi^{2}\big)^{(i+j)}|\widehat{u}_0(\xi)_{k\ell}|^2 +|\widehat{u}_1(\xi)_{k\ell}|^2 \Big) \notag \\ &  =  \mathrm{e}^{\left(- b+\sqrt{b^2-4m^2}\right) t}  \Big( \| u_0\|^2_{H^{i+j}_\mathcal{L}(\mathbb{G)}} +\| u_1\|^2_{L^2(\mathbb{G)}}\Big) \label{proof L^2 est u(t) case b^2>4m^2}
\end{align} for any $i,j\in\{0,1\}$ such that $0\leqslant i+j\leqslant 1$ and any $t\geqslant 0$. It is easy to see that \eqref{proof L^2 est u(t) case b^2>4m^2} implies \eqref{L^2 norm u(t)}, \eqref{L^2 norm (-L)^1/2 u(t)} and \eqref{L^2 norm ut(t)} in the case $b^2>4m^2$.

\section{Proof of Theorem \ref{Thm glob esistence}} \label{Subsection Fixed Point}

A fundamental tool to prove the global existence result is the following Gagliardo -- Nirenberg type inequality, whose proof can be found in \cite{RY19} (see also \cite[Corollary 2.3]{Pal20WEd}).

\begin{lemma}\label{Lemma GN ineq L2}
Let $\mathbb{G}$ be a connected unimodular Lie group with topological dimension $n\geqslant 3$. For any $q\geqslant 2$ such that $q\leqslant \frac{2n}{n-2}$ the following Gagliardo -- Nirenberg type inequality holds
\begin{align} \label{GN ineq L2}
\| f\|_{L^q(\mathbb{G})}\lesssim \| f\|^{\theta(n,q)}_{H^{1}_\mathcal{L}(\mathbb{G})} \| f\|^{1-\theta(n,q)}_{L^{2}(\mathbb{G})}
\end{align}
for any $f\in H^{1}_\mathcal{L}(\mathbb{G})$, where $ \theta(n,q)\doteq  n\left(\frac{1}{2}-\frac{1}{q}\right)$.
\end{lemma}

As we explained in Section \ref{Section approach}, in order to prove Theorem \ref{Thm glob esistence} it suffices to show the validity of \eqref{1st ineq N op} and \eqref{2nd ineq N op}.

It is helpful to rewrite $Nu=  u^{\mathrm{ln}}+ J u$, where 
\begin{align*} 
u^{\mathrm{ln}}(t,x) & \doteq  u_0(x)\ast_{(x)} E_0(t,x;b,m^2) + u_1(x)\ast_{(x)} E_1(t,x;b,m^2), \\ J u(t,x) & \doteq \int_0^t |u(s,x)|^p\ast_{(x)} E_1(t-s,x;b,m^2) \, \mathrm{d}s.
\end{align*} Let us get started by estimating $\|Nu\|_{X(T)}$ for $u\in X(T)$. By Proposition \ref{Prop L^2-L^2 estimates} it results 
\begin{align}\label{estimate u ln in X(T)}
\|u^{\mathrm{ln}}\|_{X(T)}\lesssim  \|(u_0,u_1)\|_{H^1_\mathcal{L}(\mathbb{G})\times L^2(\mathbb{G})}.
\end{align} On the other hand, thanks to the invariance by time translations of \eqref{linear wave CP}, combining Proposition \ref{Prop L^2-L^2 estimates}, the Gagliardo -- Nirenberg type inequality from Lemma \ref{Lemma GN ineq L2} and \eqref{def norm X(T)}, we get
\begin{align}
\|\partial_t^j (-\mathcal{L})^{i/2} Ju(t,\cdot)\|_{L^2(\mathbb{G})} & \lesssim \int_0^t d_{b,m^2}(t-s) \| u(s,\cdot)\|^p_{L^{2p}(\mathbb{G})} \, \mathrm{d}s \notag \\ 
& \lesssim \int_0^t d_{b,m^2}(t-s) \| u(s,\cdot)\|^{p\theta(n,2p)}_{H^{1}_{\mathcal{L}}(\mathbb{G})} \| u(s,\cdot)\|^{p(1-\theta(n,2p))}_{L^2(\mathbb{G})} \, \mathrm{d}s \notag \\ 
& \lesssim  \int_0^t d_{b,m^2}(t-s) (d_{b,m^2}(s))^p  \, \mathrm{d}s \, \| u\|^{p}_{X(t)} \notag \\ &\lesssim 
d_{b,m^2}(t)  \, \| u\|_{X(t)}^p \label{estimate Ju in X(T)}
\end{align} for $i,j\in\{0,1\}$ such that $0\leqslant i+j\leqslant 1$. We underline that the employment of \eqref{GN ineq L2} in the second step of the previous chain of inequality, we have to require the condition $p\leqslant \frac{n}{n-2}$ in Theorem \ref{Thm glob esistence}. 
So, from \eqref{estimate u ln in X(T)} and \eqref{estimate Ju in X(T)} we obtain \eqref{1st ineq N op}.

Let us derive now \eqref{2nd ineq N op}. Combining H\"older's inequality and the inequality $$||u|^p-|v|^p|\leqslant p |u-v|(|u|^{p-1}+|v|^{p-1}),$$ we get 
\begin{align*}
\| |u(s,\cdot)|^p-|v(s,\cdot)|^p\|_{L^{2}(\mathbb{G})} \lesssim \| u(s,\cdot)-v(s,\cdot)\|_{L^{2p}(\mathbb{G})} \left(\|u(s,\cdot)\|^{p-1}_{L^{2p}(\mathbb{G})}+ \|v(s,\cdot)\|^{p-1}_{L^{2p}(\mathbb{G})}\right).
\end{align*}
 Applying the previous inequality and \eqref{GN ineq L2} to estimate each $L^{2p}(\mathbb{G})$ -- norm that appears on the right -- hand side, we find
\begin{align}
\|\partial_t^j (-\mathcal{L})^{i/2}  & (Ju(t,\cdot)-Jv(t,\cdot))\|_{L^2(\mathbb{G})}   \notag \\ &
 \lesssim \int_0^t d_{b,m^2}(t-s)  \| |u(s,\cdot)|^p-|v(s,\cdot)|^p\|_{L^{2}(\mathbb{G})} \, \mathrm{d}s \notag \\   & \lesssim \int_0^t d_{b,m^2}(t-s)  (d_{b,m^2}(s))^p \, \mathrm{d}s \, \| u-v\|_{X(t)}\left(\|u\|^{p-1}_{X(t)}+\|v\|^{p-1}_{X(t)}\right)  \notag \\ 
& \lesssim d_{b,m^2}(t) \, \| u-v\|_{X(t)}\left(\|u\|^{p-1}_{X(t)}+\|v\|^{p-1}_{X(t)}\right) \label{estimate Ju -Jv in X(T)}
\end{align}  for $i,j\in\{0,1\}$ such that $0\leqslant i+j\leqslant 1$. From \eqref{def norm X(T)} and  \eqref{estimate Ju -Jv in X(T)} it follows \eqref{2nd ineq N op}. 

Note that thanks to the exponential decay rate $d_{b,m^2}(t)$ both in \eqref{estimate Ju in X(T)} and \eqref{estimate Ju -Jv in X(T)} we have the uniform boundedness of the integral 
$$(d_{b,m^2}(t))^{-1}\int_0^t d_{b,m^2}(t-s)  (d_{b,m^2}(s))^p \, \mathrm{d}s$$
without requiring further conditions on $p$.

\section{Final remarks} \label{Section final rem}

In \cite{Pal20WEd,Pal20WE} it is emphasized how the global dimension of $\mathbb{G}$ (which is 0 in the compact case) has a relevant role in the search for globally defined solution for the semilinear Cauchy problem with power nonlinearity associated to the damped wave operator  $\partial_t^2-\mathcal{L}+\partial_t$ and to the wave operator  $\partial_t^2-\mathcal{L}$, respectively. In both cases, we may not prove any global existence result. On the contrary, for any $p>1$ local in time solutions to these semilinear problem blow up in finite time under suitable sign assumption for the Cauchy data. In this paper, we showed that the contemporary presence of a damping term and of a mass term reverse completely the situation. 

As we mentioned in the introduction, this fact has been already observed in the Euclidean case (also in the case with fractional in space operators 	\cite{DAb15} or with suitable time -- dependent coefficients for the damping and mass terms \cite{Gir19}) and for graded Lie group. Nevertheless, it is remarkable to observe this phenomenon even in compact (and connected) Lie groups in consideration of what remarked above for the damped wave operator and for the wave operator, for which the corresponding semilinear Cauchy problems admit a finite critical exponent in the Euclidean case (Fujita exponent \cite{Mat76,TY01,Zhang01,IT05} and Strauss exponent \cite{John79,Glas81B,Glas81,Sid84,Scha85pn,LinSog95,Geo97,YZ06,Zhou07}, respectively) or in the special case of the damped wave equation in the Heisenberg group (see \cite{GP19DW,Pal19}).



\end{document}